\documentclass[journal,transmag, twoside]{IEEEtran}
\usepackage{float}

\usepackage{cite}

\ifCLASSINFOpdf
  \usepackage[pdftex]{graphicx}
\else
\fi
\usepackage{amsmath}
\usepackage{amssymb}
\usepackage{array}
\usepackage{fixltx2e}

\usepackage{stfloats}
\usepackage{url}

\usepackage{caption}
\usepackage{mathtools}
\usepackage{placeins}

\begin{document}
%
\title{Mathematical Analysis of Modified BEM-FEM Coupling Approach for 3D Electromagnetic Levitation Problem}


\author{\IEEEauthorblockN{Sayan Sarkar\IEEEauthorrefmark{1},
Amit Jena\IEEEauthorrefmark{2},
and Mashuq-un-Nabi\IEEEauthorrefmark{3}}
\IEEEauthorblockA{\IEEEauthorrefmark{1}Electronics and Computer  Engineering Department, Hong Kong University of Science and Technology, Hong Kong}
\IEEEauthorblockA{\IEEEauthorrefmark{2} Electrical and Computer Engineering Department , Texas A\&M University, USA}
\IEEEauthorblockA{\IEEEauthorrefmark{3}Electrical Engineering Department, IIT Delhi, New Delhi, India}
\thanks{Manuscript  received  May  x,  xxxx;  revised  August  xx,  xxxx;  accepted August xx, xxxx. Date of current version November xx, xxxx.Corresponding author: Sayan Sarkar (email:ssarkar@connect.ust.hk). \newline Color  versions  of  one  or  more  of  the  figures  in  this  article  are  availableonline at http://ieeexplore.ieee.org.
\newline Digital  Object Identifier  xx.xxxx/TMAG.xxxx.xxxxxxx}}


\markboth{IEEE TRANSACTIONS ON MAGNETICS, VOL. xx, NO. xx, DECEMBER xxxx}{{SARKAR} \MakeLowercase{et al.}:Mathematical Analysis of Modified BEM-FEM Coupling Approach for 3D Electromagnetic Levitation Problem}

%


\IEEEpubid{\begin{minipage}{\textwidth}\ \\[12pt]\centering
  xxxx--xxxx \copyright~2019 IEEE. Personal use is permitted, but republication/redistribution requires IEEE permission.\\ 
 See http://www.ieee.org/publications standards/publications/rights/index.html for more information
\end{minipage}}


\IEEEtitleabstractindextext{%
\begin{abstract}

In electromagnetic analysis, the finite element and boundary element methods jointly known as 'FEM-BEM coupling' is applied for numerically solving levitation problem based on eddy current. The main focus behind this coupled analysis method is to determine the dynamic characteristic of the levitating body in the presence of a magnetic field. An innovative 3D structure is developed that couples Lagrangian description and BEM-FEM coupling method for this purpose. The coupling methodology is based on the boundary conditions on the common boundaries between FEM and BEM sub-domains. Subsequent coding has been developed to simulate the problem in the MATLAB environment. An example similar to TEAM (Testing Electromagnetic Analysis Methods) workshop problem 28 has been used to study the efficiency of code for computationally inexpensive analysis.

\end{abstract}

}

\maketitle

\IEEEdisplaynontitleabstractindextext

%
\IEEEpeerreviewmaketitle

\section{Introduction}
%
%
%
%

\IEEEPARstart{T}{he} modelling of dynamics of a conducting object during external magnetic field's presence can be calculated by solving the equation of motion considering the magnetic forces for induced eddy currents. Traditional difficulties in this area are the analysis of the motion of conducting body, re-meshing techniques, force calculation, weak versus strong electromechanical coupling, efficient time-stepping schemes etc. \cite{karl},\cite{kurz1}. TEAM (Testing Electromagnetics Analysis Method) workshop problems 9 and 17 manage moving bodies with external magnetic field\cite{jena}. In contrast, TEAM workshop problem 28 is a transient problem with electromechanical coupling. The proposed model for evaluating levitating setup is similar to TEAM workshop problem 28. 

In general, the discretization of a domain by finite elements which contain movable parts is not straightforward to handle\cite{kurz1}. Even a small displacement can create huge mesh distortion; subsequent re-meshing of the entire domain is necessary for tackling the continuous displacement of the moving parts. This shortcoming could be removed \cite{fetter} by the BEM-FEM coupling approach. The advantage of this approach is meshing pattern does not alter during the motion. It gets shifted according to the levitating object due to BEM application and there is no need for entire domain re-meshing. The behavior of the levitated body is described by Maxwell's equations and other constitutive equations. Solving the Electromechanical problem by those equations in three dimensional space is a time-consuming CPU task, consequently barely resorted \cite{suuriniemi}. The ongoing increment of computation capability and constantly changing numerical techniques have dramatically diminished the need for time-consuming calculations. During evaluation, 'BEM-FEM' coupling code has been developed in MATLAB environment. Laplace equation is used as governing equation in the air (BEM) domain for simulation. Results have been compared with measured values from Multiphysics software "COMSOL."An analysis without detailed mathematical treatment of the BEM approach is provided in \cite{jena} to calculate the magnetic vector potential at the boundary of the conducting medium. In this paper, an extension of previous work\cite{jena}, based on mathematical ground is provided.

\section{ELECTROMAGNETIC FORMULATION}
Maxwell’s equations that describes all classical electromagnetic phenomena, are given by (1)-(4).
\begin{equation}
\nabla \cdot \vec{E}=\frac \rho {\epsilon}
\label{eqn:nab1}
\end{equation}
\begin{equation}
\nabla \cdot \vec{B}=0
\label{eqn:nab2}
\end{equation}
\begin{equation}
\nabla \times \vec{E}
	= - \frac {\partial{\vec B}} {\partial t}
\label{eqn:nab3}
\end{equation}
\begin{equation}
\nabla \times \vec{B}
	= \mu(\vec J+\epsilon\frac {\partial{\vec E}} {\partial t})
\label{eqn:nab4}
\end{equation}
Displacement current is omitted for sake of simplifications.Above (\ref{eqn:nab3}) and (\ref{eqn:nab4}) displayed an crucial role for development of transient electromagnetic problem.Useful identity of vector calculus gives (5)
\begin{equation}
\vec{E}= -\frac{\partial{\vec A}} {\partial t}-\nabla \vec{V}
\label{eqn:vec1}
\end{equation}
\emph'{A}' and \emph '{V}' are magnetic vector potential and electric scalar potential respectively. Object under study has no residual magnetism. Accounting constitutive relations in (\ref{eqn:vec2}) and (\ref{eqn:vec3}).
\begin{equation}
\vec{B}=\mu\vec{H}
\label{eqn:vec2}
\end{equation}
\pagebreak
\begin{equation}
\vec{J}=\vec{J_s}+\sigma(\vec{E}+\vec{\nu}\times\vec{B})
\label{eqn:vec3}
\end{equation}

Here, $\vec{J_s}$, $\vec{H}$, $\sigma$ and $\vec{\nu}$ are impressed current density, magnetic field strength,  conductivity and velocity respectively. Substituting (\ref{eqn:vec1}) in (\ref{eqn:vec3}) gives,

\begin{equation}
\vec{J}-\sigma(-\frac{\partial{\vec{A}}}{\partial t}-\nabla\vec{V}) - \sigma(\vec{\nu}\times\vec{B})=\vec{J_s}
\label{eqn:vec4}
\end{equation}
Simplification of (\ref{eqn:vec4}) yields

\begin{equation}
-\frac{1}{\mu}\nabla^2 \vec{A}+\sigma{\frac{\partial \vec A}{\partial t}+\nabla \vec V-(\vec \nu\times \vec B)}=\vec{J_s}
\label{eqn:11}
\end{equation}
The divergence of current density is also zero.
\begin{equation}
\nabla\cdot \vec J=0
\label{eqn:12}
\end{equation}
Substituting (10) into (8) gives (11).
\begin{equation}
\nabla\cdot\sigma{\frac{\partial\vec A}{\partial t}+\nabla\cdot \vec V}-(\vec\nu\times \vec B)=0
\label{eqn:13}
\end{equation}

\begin{figure}[h]
\includegraphics[width=8cm]{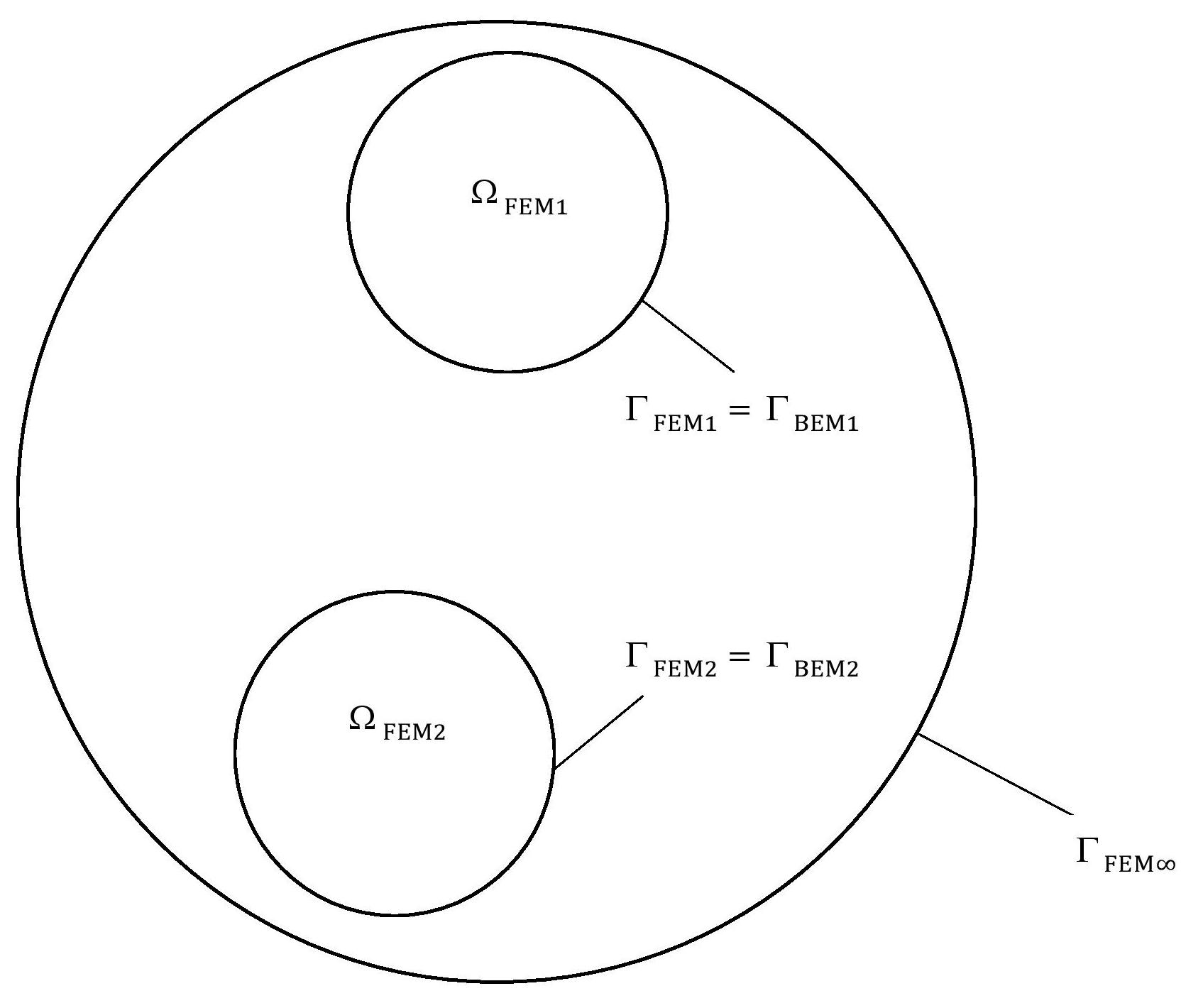}
\centering
\caption{structure of the considered domain\cite{jena}}
\label{fig:bf}
\end{figure}

The considered domain $\Omega$ is decomposed as shown in Fig. \ref{fig:bf} into different conducting region like $\Omega_{FEM1}$, a region $\Omega_{FEM2}$ free of eddy current (air domain). In multiply connected region $\Omega_{BEM}$, there is no conducting body but it may contain a current source, The BEM treats the surrounding air space. The following continuity conditions are valid on coupling boundaries for magnetic vector potential, magnetic field intensity. Current density is continuous across the coupling boundary.

\begin{equation}
\vec{A}_{FEM}=\vec{A}_{BEM}
\label{eqn:14}
\end{equation}
\begin{equation}
\nabla\cdot \vec{A}=0\quad \text{(Coulomb\, gauge\, condition)}
\label{eqn:14a}
\end{equation}
\begin{equation}
\nabla\cdot\vec A_{FEM}=\nabla\cdot\vec A_{BEM}
\label{eqn:15}
\end{equation}
\begin{equation}
\frac{1}{\mu}\nabla\times(\vec A\times \vec n)\vert^{FEM}=\frac{1}{\mu}\nabla\times(\vec A\times \vec n)\vert^{BEM}
\end{equation}

\section{ELECTROMECHANICAL SYSTEM MODELING}

\begin{figure}[H]
\includegraphics[width=6cm]{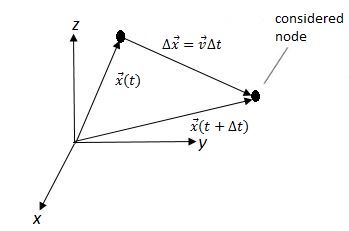}
\centering
\caption{Node moving from $\Delta$\emph{x} to \emph{x}+$\Delta \emph{x}$ within $\Delta$\emph{t}\cite{jena}}
\label{fig:vd}
\end{figure}

In the general case of movable boundaries, at least a portion of the meshing structure must move as a whole \cite{kurz2} with the boundary to account for the interface conditions. Fig. \ref{fig:vd} shows a node which is part of a moving mesh. The temporal variation of the vector potential at this node is mentioned by \cite{fetter}.

\begin{equation}
\frac{d\vec A}{dt}=\frac{\partial \vec A}{\partial t}+(\vec\nu\cdot\nabla)\vec A
\label{eqn:17}
\end{equation}

Assuming nodal elements, the temporal variation of the magnetic vector potential values gives an approximation for the total time derivative rather than the partial time derivative \cite{fetter}. Hence
\begin{equation}
(\vec\nu\times\vec B)=\nabla(\vec\nu\cdot A)-\frac{d\vec A}{dt}+\frac{\partial\vec A}{\partial t}-\vec A\times \vec\omega
\label{eqn:18}
\end{equation}
substituting (\ref{eqn:18}) in (\ref{eqn:11}), we have
\begin{equation}
-\frac{1}{\mu}\nabla^2 \vec A+\sigma\left\{ {\nabla( V^*)+\frac{d\vec A}{dt}} +\vec A\times \vec\omega\right\}=\vec J_s
\label{eqn:19}
\end{equation}
similarly substituting (\ref{eqn:18}) in (\ref{eqn:13}) 
\begin{equation}
\nabla\cdot\sigma\left\{ \nabla(V^*)+\frac{d\vec A}{dt}+\vec A\times\vec\omega\right\}=0
\label{eqn:20}
\end{equation}
\begin{center}
Where $\vec V^*$=$V-\nu A$
\end{center}

From (\ref{eqn:19})-(\ref{eqn:20}), it gives that moving the mesh \cite{fetzer} and replacing \emph{V} by \emph{$V^*$} automatically takes into account the subsequent motional E.M.F till $\vec A\times\vec\omega=0$\cite{fetter}. This is a scenario for 2-D problems and 3-D problems with small perturbation. In those cases, it is completely sufficient to move the mesh connected to the respective moving bodies for all electromagnetic effects due to motion.

\section{FEM-BEM COUPLING METHOD}

The application of the Galerkin method for the weak integral form of (\ref{eqn:19})-(\ref{eqn:20}) and using nodal elements and boundary conditions yields the coupled matrix of FEM-BEM method.FEM-BEM method across coupling boundaries is given by (\ref{eqn:21}).

\begin{equation}
\begin{multlined}
[\frac{1}{\mu_1}(\nabla\times \vec A)\times\vec n_1+\frac{1}{\mu_0}(\nabla\times \vec A)\times\vec n_2\\
+\frac{1}{\mu_1}(\nabla\cdot\vec A)\vert_1 n_1+\frac{1}{\mu_0}(\nabla\cdot\vec A)\vert_2 n_2 ]=0
\end{multlined}
\label{eqn:21}
\end{equation}

In order to apply Galerkin in (\ref{eqn:19}) and (\ref{eqn:20}), we multiply (18) with a suitable weight function.

$\omega_a$
\begin{equation}
\begin{multlined}
\int\limits_{\Omega}[\frac{1}{\mu}\nabla\times(\nabla\times\vec  A)-\nabla\cdot(\frac{1}{\mu_0}\nabla\cdot\vec A)\\
+\sigma\left\{\frac{d\vec A}{dt}+\nabla \vec V^* +\vec A\times\vec \omega \right\}]\omega_ad\Omega=\int\limits_{\Omega}J_s\omega_{a}d\Omega
\end{multlined}
\label{eqn:22}
\end{equation}
\begin{center}
where, a = 1, 2, 3.
\end{center}
\begin{equation}
\int\limits_{\Omega}\nabla\cdot\sigma\left\{\frac{d\vec A}{dt}+\nabla V^*+\vec A\times\vec \omega\right\}\omega d\Omega=0
\label{eqn:23}
\end{equation}
Where,
\begin{center}
$
\omega_1=
\begin{pmatrix}
\omega_1\\0\\0
\end{pmatrix},
\omega_2=
\begin{pmatrix}
0\\\omega_2\\0 
\end{pmatrix},
\omega_3=
\begin{pmatrix}
0\\0\\\omega_3
\end{pmatrix}
$
\end{center}

Application of boundary condition,vector identity, weight function and further simplification yields (\ref{eqn:24a}).

\begin{equation}
\begin{multlined}
\int\limits_{\Omega_1}[\nabla\cdot\sigma\left\{\frac{\partial\vec A}{\partial t}+\nabla V^*+\vec\omega\times \vec A\right\}]\\
=-\int\limits_{\Omega_1}\sigma\left\{\frac{\partial \vec A}{\partial t}+\nabla V^* +\vec\omega\times \vec A\right\}\nabla\omega d\Omega\\
+\int\limits_{\Gamma_{12}}\sigma\left\{\nabla(V^*)+\frac{d\vec A}{dt}+\vec A\times\vec\omega\right\}n_1\omega d\Gamma
\end{multlined}
\label{eqn:24a}
\end{equation}

Simplification gives (24).

\begin{equation}
\int\limits_{\Omega}\nabla\cdot\sigma\left\{\frac{\partial\vec A}{\partial t}+\nabla V^*+\vec\omega\times\vec A\right\}\omega d\Omega = 0
\label{eqn:25}
\end{equation}

It can be further simplified again using vector analysis, and taken into account $\sigma=0$ and $M=0$ and $\mu=\mu_0$ it can be written as follows-

\begin{equation}
\begin{multlined}
-\frac{1}{\mu_0}\int_{\Omega_{BEM}}(\nabla^2\omega_a)Ad\Omega-\int_{\Gamma_{BEM}}\frac{1}{\mu_0}\left\{(\vec n_2\times\nabla)\times\omega_a\right.\\
\left.-(\vec n_2\times\nabla)\times\omega_a\right\}Ad\Gamma-\int\limits_{\Gamma_{BEM}}\vec Q\cdot\omega_ad\Gamma\\
=\int\limits_{\Omega_{BEM}}J_s\cdot\omega_a\cdot d\Omega
\end{multlined}
\label{eqn:27}
\end{equation}
The fundamental solution of scalar Laplace equation in 3-D is given by following equation -
\begin{equation}
\Delta u^*=-\delta(r-r')
\label{eqn:28}
\end{equation}

It gives therefore  
\begin{equation}
\omega_1=
\begin{pmatrix}
u^*\\0\\0 
\end{pmatrix},
\omega_2=
\begin{pmatrix}
0\\u^*\\0 
\end{pmatrix},
\omega_3=
\begin{pmatrix}
0\\0\\u^*
\end{pmatrix}
\label{eqn:29}
\end{equation}
By substituting (\ref{eqn:29}) in (\ref{eqn:27}) gives (\ref{eqn:30}).
\begin{equation}.
\begin{multlined}
\frac{1}{\mu_0}c(r)A+\frac{1}{\mu_0}\int\limits_{\Gamma_{BEM}}((\vec n_2\times\nabla u^*)\times\vec A+q^*A)
d\Gamma\\-\int\limits_{\Gamma_{BEM}}u^*Qd\Gamma=\int\limits_{\Omega_{BEM}}J_s u^*d\Omega
\end{multlined}
\label{eqn:30}
\end{equation}

Where, $q^*=\frac{\partial u}{\partial n_2}$; $c(r)=\frac{\theta}{4\pi}$ with $c(r)$ being called as edge factor and $\theta$  is the solid angle at the point rim area inside. Quantity c(r) is equal to 1 if the field point belongs to the interior point. If it lies on the smooth part of the boundary, then it is equal to $\frac{1}{2}$, and it takes 0 value if the point lies outside. Till now, the boundary element method-based integral equation has been derived whose solution can be obtained by clustering of the boundary $\Gamma _{BEM}$ and selection of appropriate basis function A and $Q_1$ \cite{kurz3}. The edge $\Gamma_{BEM}$ is in several non-core elements $\Gamma_e$ decomposed.The basic condition for coupling along boundary \cite{fetter}- \cite{kurz2}.

\begin{equation}
\begin{cases}
\vec H_t\vert_{FEM}&=\vec H_t\vert_{BEM}\\
\nabla\cdot\vec A_{FEM}&=\nabla\cdot\vec A_{BEM}\\
\vec A_{FEM}&=\vec A_{BEM}\\
\vec B_n\vert_{FEM}&=\vec B_n\vert_{BEM}\\
\vec\Gamma_{BEM}&=\vec\Gamma_{FEM}
\end{cases}
\label{eqn:31}
\end{equation}
\begin{equation}
\begin{multlined}
\frac{1}{\mu_0}c(r_k)A+\frac{1}{\mu_0}\int\limits_{\Gamma_{BEM}}q^*Ad\Gamma\\
-\sum_{e=1}^{E}\int\limits_{\Gamma_{BEM}}u^*Q_1d\Gamma
=\int\limits_{\Omega_{BEM}}J_su^*d\Omega
\end{multlined}
\label{eqn:32}
\end{equation}

The vector element A and $Q_1$ interpolated to $\Gamma_e$ (boundary element) to evaluate the integrals (\ref{eqn:32}) numerically \cite{jena}. The functions used here are the same as in the finite element method. The standard Lagrangian method and subsequent transformation reduced the dimension of the problem by one unit. The matrix format of (30) follows (31).

\begin{equation}
[H]\left\{A\right\}-[G]\left\{Q_1\right\}=\left\{F(J)\right\}
\label{eqn:34}
\end{equation}
Similarly the FEM domain written in matrix form gives (32).
\begin{equation}
[k]\left\{A\right\}-[t][Q_1]=\left\{F(J)\right\}
\label{eqn:35}
\end{equation}
Where in domain $\Omega_e$
\begin{center}
$A_e(r,t)$\qquad=\:$\sum_{j=1}^{n}N_j(r)A_j(t)$\\[1ex]
$Q_{1e}(r,t)$\qquad=\:$\sum_{j=1}^{n}N_j(r)Q_{1e}(t)$\\[1ex]
$V_e(r,t)$\qquad=\:$\sum_{j=1}^n N_j(r)V_j(t)$
\end{center}
Now combination of (\ref{eqn:34}) and (\ref{eqn:35}) with above condition yields the local matrix of FEM-BEM coupling.
\begin{equation}
[c]\left\{A, V^*\right\}+[k]\left\{A, V^*\right\}-[t]\left\{Q_1\right\}=\left\{F(J)\right\}
\label{eqn:36}
\end{equation}
Where k, c, t are respectively element stiffness matrix, damping matrix and boundary matrix.The element equation for $e=1, 2, 3, ..., L$ are obtained for obtaining the global matrix for FEM-BEM coupling \cite{kurz3}
\begin{equation}
\begin{multlined}
[c]\left\{A,V^*\right\}+\left\{[k]+([T][G]^{-1}[H])\right\}\left\{A, V^*\right\}\\
=[T][G]^{-1}\left\{RS\right\}+\left\{F(J)\right\}
\end{multlined}
\label{eqn:41}
\end{equation}
Now, the quantity under first and second bracket in (\ref{eqn:41}) is denoted by $K _{BEM}$ and $K_{TOTAL}$ respectively where
$K_{BEM}$ refers to matrix equivalent of finite element and  (\ref{eqn:41}) is a linear differential equation of first order\cite{kurz3} . 

\section{MODIFIED APPROACH OF FEM-BEM COUPLING}
The traditional approach of BEM-FEM coupling in 3D proceeds as stated in the previous section. The simultaneous implementation of the FEM and BEM approach raises the computational load and execution time. According to the modified approach proposed here, FEM has to be implemented once during the whole process in the area encircled by the internal circle as shown in Fig. \ref{fig:set1}-\ref{fig:set2}. At external circle in Fig. \ref{fig:set1}-\ref{fig:set2}, boundary condition are nearly zero as stated earlier. FEM provides as the input of BEM and due to lack of magnetization vector in Aluminium body. The output of FEM remains almost constant during the levitation process. FEM has been adopted once in the internal circle domain, and output is obtained on the boundary of the internal circle as in Fig.\ref{fig:set1}-\ref{fig:set2}. BEM is used as the propagation tool to obtain the required result at the boundary of the metallic object. In order to make the calculations simpler, the entire 3D object has been projected to 2D, and analysis, computations have been done in 2D. It is analogous to the realization that the planar view of the levitation process is observed and the depth of the body has been taken as unity. In planar view, the governing equation for FEM and BEM sub-domains is converted to 2D Laplace and 2D poison's equation, respectively. 
\newline In mathematics for differential equations, a boundary value problem is a differential equation together with a set of additional constraints, called the boundary conditions. A solution to a boundary value problem is a solution to the differential equation, satisfying the boundary conditions. Dirichlet (or first-type) boundary condition is a type of boundary condition when imposed on an ordinary or a partial differential equation, it gives those values which need to take on along the boundary of the domain. Neumann (or second-type) boundary condition is a type of boundary condition imposed on an ordinary or partial differential equation. It specifies derivative based values of a solution on the domain's boundary. Robin boundary conditions are a weighted combination of Dirichlet boundary conditions and Neumann boundary conditions (mixed boundary condition). In contrast to mixed boundary conditions, only one type of boundary conditions are applied on any portion of the boundary traditionally: either the function value is specified or the normal derivative but not both. 

\begin{equation}
xA+y\frac{\partial A}{\partial n}=g\qquad on \quad \partial\Omega
\label{41a}
\end{equation}
For some non-zero constants, x and y and a predefined function $g$ defined on $\partial\Omega$. Here, $A$ is the unknown solution defined on $\Omega$, and $\frac{\partial A}{\partial n}$ denotes the normal derivative at the boundary. More generally, x and y are functions rather than constants.In the BEM sub-domain the governing equation is given by \cite{ang}.

\begin{equation}
\nabla^2 A =0
\label{eqn:42}
\end{equation}

Fundamental solution of '$A$' plays an pivotal role in formation of boundary integral equation and it is given by (37).Where (x, y) is the source point and $(\xi, \eta)$ is the field point.
\begin{equation}
A(x,y,\xi,\eta) = \frac{ln((x-\xi)^2+(y-\eta)^2)}{4\pi}
\label{field}
\end{equation}
 From Gauss Divergence theorem and (\ref{eqn:42}) it can be stated that
 \begin{equation}
\int\limits_{C}(A_2 \frac{\partial A_1}{\partial n}-A_1\frac{\partial A_2}{\partial n})ds(x,y)=0
\label{eqn:43}
\end{equation}
$A_1$ and $A_2$ are two arbitrary solutions of Laplace's equation in region R bounded by simple closed curved C.$\frac{\partial A_1}{\partial n}$ and $\frac{\partial A_2}{\partial n}$ are normal derivatives of the functions. Fundamental solution and the general solution satisfy (\ref{eqn:42}).Hence,

\begin{equation}
\scalebox{0.80}{$
\begin{aligned}
\int\limits_{C}[A(x,y)]\frac{\partial}{\partial n}A(x,y;\xi,\eta)-A(x,y;\xi,\eta)\frac{\partial}{\partial n}A(x,y)]ds(x,y) = 0
\end{aligned}$}
\label{eqn:44}
\end{equation}

In (\ref{eqn:44}), the fundamental solution possesses a discontinuity at $x=\xi$, $y=\eta$. Hence, $(\xi, \eta) \notin R\cup C$ 
\begin{figure}[h]
\includegraphics[width=6cm]{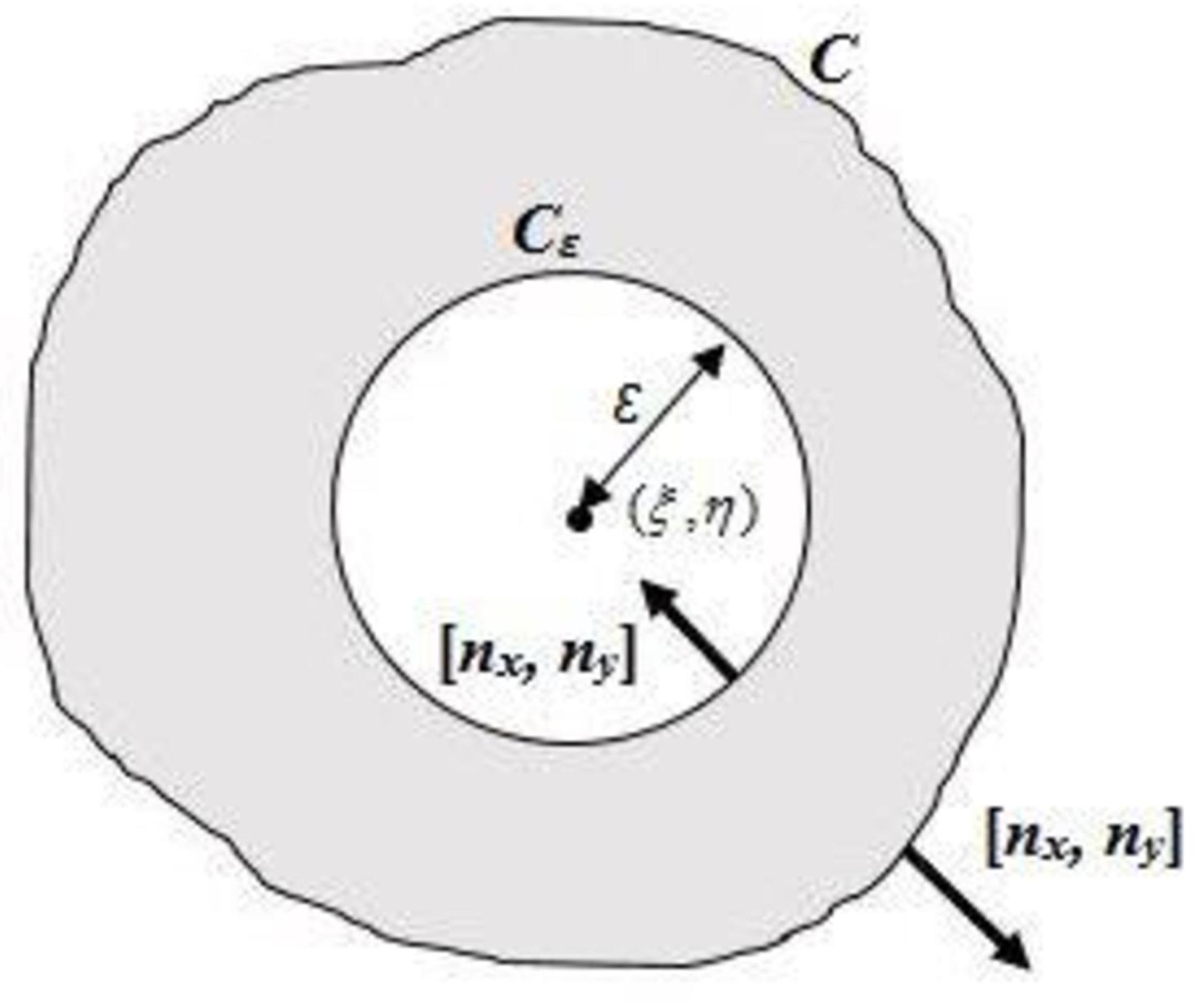}
\centering
\caption{Modified approach when a field lie in the interior\cite{jena}}
\label{fig:ab1}
\end{figure}
If $(\xi, \eta)$ $\in$ R then C has to be replaced by C $\cup$ $C_\epsilon$ \cite{ang}.$C_\epsilon$ is a circle of center $(\xi, \eta)$ radius $\epsilon$ has center as shown in Fig. \ref{fig:ab1}. This is because $\phi(x,y,,\xi,\eta)$ along with its partial derivatives are well defined in the annular region between $C$ and $C_\epsilon$.

\begin{figure}[H]
\includegraphics[width=4.5cm]{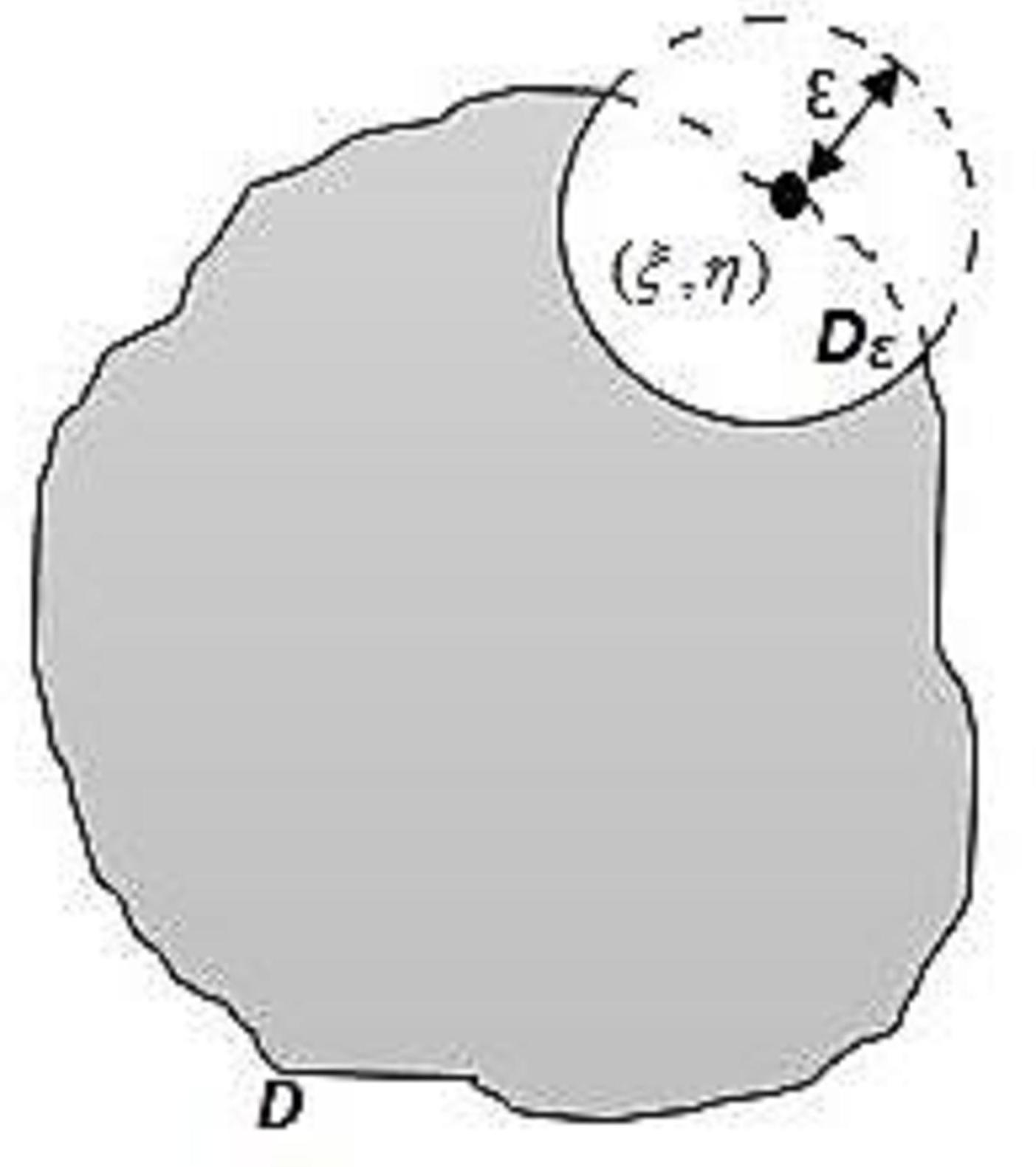}
\centering
\caption{Modified approach when a field lie on the boundary \cite{jena}}
\label{fig:ab2}
\end{figure}
If $(\xi, \eta)$ $\in$ C lies on the boundary, then C has to be changed by $D\cup D_\epsilon$ as in Fig. \ref{fig:ab2}. If $C_\epsilon$ has the center $(\xi, \eta)$ and radius $\epsilon$, then D is the part of C that lies outside $C_\epsilon$.
$D_\epsilon$ is that of $C_\epsilon$ that is inside R. Actually to solve the equations numerically, we need to discretize the boundary integral equation. That is

\begin{equation}
\begin{multlined}
\lambda(\xi,\eta)A(\xi,\eta)=\int\limits_{C}[A(x,y)\frac{\partial}{\partial n}A(x,y; \xi,\eta)\\
-A(x,y; \xi,\eta)\frac{\partial}{\partial n}A(x,y)]ds(x,y)
\end{multlined}
\label{eqn:45}
\end{equation}
Where $(\xi,\eta)$ is any arbitrary point in $R\cup C$ and
\begin{align*}
\lambda(\xi,\eta)=
\begin{cases}
0 &\qquad if \: (\xi,\eta)\notin R\cup C\\
\frac{1}{2}&\qquad \text{if}\: (\xi,\eta)\: \text{on}\: \text{the}\:\text {smooth}\:\text{part}\: \text{of}\: C\\
1&\qquad \text{if} \: (\xi,\eta)\in R
\end{cases}
\end{align*}
In this approximation procedure of $A$ and $\frac{\partial A}{\partial n}$ on the boundary condition of C is required. Here we have taken constant approximation, i.e., it is assumed that $A$ and $\frac{\partial A}{\partial n }$ functions are constants over $C_K$ where $C\approx C_1\cup C_2\cup C_3\cup C_4...\cup C_{N-1}\cup C_N$. More precisely it can be said that $A=\bar{A}_K$ and $\frac{\partial A}{\partial n}=\frac{\partial \bar{A}_K}{\partial n}$ where $\bar{A}_K$ and $\frac{\partial \bar{A}_K}{\partial n}$ are values of $A$ and $\frac{\partial A}{\partial n}$ at midpoints of $C_K$. Here, $C_K$ is one arc of the boundary or simply a boundary element. Now, the boundary integral equation boils down to (\ref{eqn:46}).
\begin{equation}
\begin{multlined}
\lambda(\xi,\eta)A(\xi,\eta)=\bar{A}_K\int\limits_{C}\frac{\partial}{\partial n}A(x,y;\xi,\eta)ds(x,y)\\[0.5ex]
-\frac{\partial\bar{A}^K}{\partial n}\int\limits_{C}A(x,y;\xi,\eta)ds(x,y)
\end{multlined}
\label{eqn:46}
\end{equation}

Here, the interior region 'R' is the annulus created by two circles. $C$ is the boundary of a doubly connected annular domain. Hence, discretizing $C$ into 80 boundary segments where precisely 40 of those lie on the boundary of the smaller circle and the other 40 elements lie on the other one. Constant shape functions have been taken to simplify the formulation process in (\ref{eqn:46}) and it leads to (\ref{eqn:47}).

\begin{equation}
\lambda(\xi,\eta)A(\xi,\eta)=\sum_{K=1}^{80}\bar{A}_K F_2^K(\xi,\eta)-\sum_{K=1}^{80}\bar{P}_K F_1^K(\xi,\eta)
\label{eqn:47}
\end{equation}
Where,
\begin{center}
$F_1^K(\xi,\eta)=\int\limits_{C_K}A(x,y;\xi,\eta)ds(x,y)$\\[0.5ex]
$F_2^K(\xi,\eta)=\int\limits_{C_K}\frac{\partial}{\partial n}A(x,y;\xi,\eta)ds(x,y)$\\[0.5ex]
$\bar{P}_K=\frac{\partial \bar{A}_K}{\partial n}$
\end{center}
As it can be seen from (\ref{eqn:47}) both $\bar{A}_K$ and $\frac{\partial \bar{A}_K}{\partial n}$ are needed on every individual boundary element $C_K$. Here $\bar{A}_K$ is the Dirichlet boundary condition and $\frac{\partial \bar{A}_K}{\partial n}$ is the Neumann boundary condition. Generally, both are not given for any node; hence, to apply (\ref{eqn:47}), we have to find missing Dirichlet conditions where only Neumann conditions are given and vice versa. Basically, they are complementary boundary conditions. In order to get $A(\xi,\eta)$, complementary boundary conditions have to be obtained by projecting (\ref{eqn:47}) on to the boundary $C$. Hence from definition of $\lambda(\xi,\eta)$, we get (\ref{eqn:48})
\begin{equation}
\frac{1}{2}\bar{A}_m=\sum_{K=1}^{80}\left\{\bar{A}_K F_2^K(\bar{x}^M,\bar{y}^M)-\bar{P}_K F_1^K(\bar{x}^M,\bar{y}^M)\right\}
\label{eqn:48}
\end{equation}
In a matrix form (\ref{eqn:48}) can be written as (\ref{eqn:49}).

\begin{equation}
\begin{multlined}
\frac{1}{2}
\begin{pmatrix}
\bar{A}_1 \\\vdots\\\bar{A}_{80}
\end{pmatrix}_{BEM}\\[1.5ex]
=
\begin{pmatrix}
F_2^1(\bar{x}^1,\bar{y}^1) & \cdots  & F_2^{80}(\bar{x}^1,\bar{y}^1)\\
\vdots & \ddots &\vdots\\
F_2^1(\bar{x}^{80},\bar{y}^{80}) & \cdots  & F_2^{80}(\bar{x}^{80},\bar{y}^{80})
\end{pmatrix}
\begin{pmatrix}
\bar{A}_1 \\\vdots\\\bar{A}_{80}
\end{pmatrix}_{BEM}\\[1.5ex]
-
\begin{pmatrix}
F_1^1(\bar{x}^1,\bar{y}^1) & \cdots  & F_1^{80}(\bar{x}^1,\bar{y}^1)\\
\vdots & \ddots &\vdots\\
F_1^1(\bar{x}^{80},\bar{y}^{80}) & \cdots  & F_1^{80}(\bar{x}^{80},\bar{y}^{80})
\end{pmatrix}
\begin{pmatrix}
\bar{P}_1 \\\vdots\\\bar{P}_{80}
\end{pmatrix}_{BEM}\\[0.1ex]
\end{multlined}
\label{eqn:49}
\end{equation}
In (\ref{eqn:49}), $(\bar{x}^m, \bar{y}^m)$ lies on the smooth part of the boundary $\lambda$ should be taken as $\frac{1}{2}$. As it has been mentioned above, both Dirichlet and Neumann conditions are needed at every node. But generally, one of the conditions is given on each node. In (\ref{eqn:47}), all $F_1$ and $F_2$ can be calculated by known formula. Hence it can be solved for N unknowns by N equations. As some node (that is for some k) $\bar{A}_K$ is unknown, and at some other node $\bar{P}_K$ is unknown. Now coming to this case, it is easy to obtain  $\bar{A}_K$ values at every node. So, the unknown quantities are $\bar{P}_K$ values (the Neumann conditions).From (\ref{eqn:49}), it can be observed that, 80 equations and 80 unknown  $\bar{P}_K$ values. The equation in (\ref{eqn:49}) can be re-written as (45).

\begin{equation}
\sum_{K=1}^{N}a^{MK}z^{K}=\sum_{K=1}^{N}b^{MK}
\label{eqn:50}
\end{equation}

For M=1, 2, 3,......, N
\begin{align*}
a^{MK}=
\begin{cases}
F_1^{K}(\bar{x}^M,\bar{y}^M)\:\:\text{if}\:\: \phi\:\: \text{is}\:\: \text{specified}\:\: \text{over}\:\: C_K\\[0.5ex]
F_2^{K}(\bar{x}^M,\bar{y}^M)-\frac{1}{2}\delta^{MK} \:\:\text{if}\:\: \frac{\partial}{\partial n}\phi\:\: \text{is}\\
\qquad\qquad\qquad\qquad\quad\:\text{specified}\:\: \text{over}\:\: C_K
\end{cases}
\end{align*}

\begin{align*}
b^{MK}=
\begin{cases}
\bar{\phi}(k)(-F_2^K(\bar{x}^M,\bar{y}^M))+\frac{1}{2}\delta^{MK}\:\:\text{if}\:\: \phi\:\: \text{is}\\
\qquad\qquad\qquad\qquad\quad\text{specified}\:\: \text{over}\:\: C_K\\[0.5ex]
\bar{P}(k)(F_1^{K}(\bar{x}^M,\bar{y}^M)\:\:\text{if}\:\: \frac{\partial}{\partial n}\phi\:\: \text{is}\\
\qquad\qquad\qquad\qquad\quad\text{specified}\:\: \text{over}\:\: C_K
\end{cases}
\end{align*}

\begin{align*}
\delta^{MK}=
\begin{cases}
0\quad \text{if}\quad M &\neq K\\
1 \quad \text{if}\quad M&=K
\end{cases}
\end{align*}

\begin{align*}
z^K=
\begin{cases}
\bar{P}(k)\:\: \text{if}\:\: \phi\:\: \text{is}\:\: \text{specified}\:\: \text{over}\:\: C_K\\[0.5ex]
\bar{\phi}(k)\:\: \text{if}\:\: \frac{\partial \phi}{\partial n}\:\:\text{is}\:\:\text{specified}\:\: \text{over}\:\: C_K
\end{cases}
\end{align*}

Once the equation in (\ref{eqn:50}) is solved for the unknown $Z^1$, $Z^2$, $Z^3$,......, $Z^N$, the values of $A$ and $\frac{\partial A}{\partial n}$ over the element $C_K$ are given by $\bar{A}_K$ and $\frac{\partial \bar{A}_K}{\partial n}$ respectively, where K=1, 2, 3,......, N. Then, $\lambda$=1 is used to compute the value of the function at any interior point. Subsequently, $F_1^K$ and $F_2^K$ values need to be computed. The points on the element $C_K$ are defined by (\ref{eqn:51})

\begin{equation}
\left.
\begin{aligned}
x&=x^K-tl^Kn_y^K\\
y&=y^K-tl^Kn_x^K
\end{aligned}
\right\}
\text{from}\: t=0 \:\text{to}\: 1
\label{eqn:51}
\end{equation}

Where $l^K$ is the length of the $C_K$ and [$n_x^K$, $n_y^K$] = $[\frac{y^{K+1}-y^K}{l^K},\frac{x^{K+1}-x^K}{l^K}]$ is the normal vector to $C_K$ pointing away from R. For (x,y) $\in$ $C_K$, it can be found that

\begin{center}
$ds(x,y)$ = $\sqrt{dx^2+dy^2}$ = $l^Kdt$
\end{center}

Using (\ref{eqn:51}), it can be stated that
\begin{center}
$((x-\xi)^2+(y-\eta)^2)$\\[1ex]
$=(x^K-\xi)^2+t^2(l^K)^2(n_y^K)^2-2(x^K-\xi)tl^Kn_y^K$\\[1ex]
$+(y^k-\eta)^2+t^2(l^K)^2(n_x^K)^2+2(y^K-\eta)tl^Kn_x^K$
\end{center}

Distance 'S' can be expressed as

\begin{equation}
S = A_Kt^2+B_K(\xi,\eta)t+E_K(\xi,\eta)
\label{eqn:53}
\end{equation}

Where
\begin{center}
$A_K=(l^K)^2$\\[1ex]
$B_K(\xi,\eta)=[-n_y^K(x^K-\eta)+n_x^K(y^K-\eta)]\times(2l^K)$\\[1ex]
$E_K(\xi,\eta)=(x^K-\xi)^2+(y^K-\eta)^2$
\end{center}
As a whole from mathematical point of view
\begin{equation}
4A_KE_K(\xi,\eta)-[B_K(\xi,\eta)]^2=0
\label{eqn:52}
\end{equation}
$F_1(\xi,\eta)$ and $F_2(\xi,\eta)$ can be written as 
\begin{equation}
F_1^K(\xi,\eta)=\frac{l^K}{2\pi}\int_0^1ln(A_Kt^2+B_K(\xi,\eta)t+E_K(\xi,\eta))dt
\label{eqn:53}
\end{equation}

\begin{equation}
F_2^K(\xi,\eta)=\frac{1}{2\pi}\int_0^1\frac{n_x^K(x^K-\xi)+n_y^K(y-\eta)}{A_Kt^2+B_K(\xi,\eta)t+E_K(\xi,\eta)}
\end{equation}
It can be mathematically proved that if $(\xi,\eta)$ lies on the boundary element on which we are integrating, then it creates a singularity in both the integrals. So two different cases have to be considered for each integral.\\[0.5ex]
$F_1^K(\xi,\eta)$ - case 1:\\[0.5ex]
In this case, $(\xi,\eta)$ doesn't lie on the boundary element on which integration is performed. Analytically this integration can be stated as in (\ref{eqn:55})
\begin{equation}
\begin{multlined}
F_1^K(\xi,\eta)=\frac{l^K}{4\pi}\left\{2[ln(l^K)-1]-\frac{B_K(\xi,\eta)}{2A_K}ln\left|\frac{E_K(\xi,\eta)}{A_K}\right|\right. \\[1ex]
\left.\left(1+\frac{B_K(\xi,\eta)}{2A_K}\right)ln\left|1+\frac{B_K(\xi,\eta)}{A_K}+\frac{E_K(\xi,\eta)}{A_K}\right|\right.\\[1ex]
\left.+\frac{\sqrt{(4A_KE_K(\xi,\eta)-[B_K(\xi,\eta)])}}{A_K}\right.\\[1ex]
\left.\times\left[tan^{-1}\left(\frac{2A_K+B_K(\xi,\eta)}{\sqrt{(4A_KE_K(\xi,\eta)-[B_K(\xi,\eta)]^2)}}\right)\right.\right.\\[1ex]
\left.\left.-tan^{-1}\left(\frac{B_K(\xi,\eta)}{\sqrt{(4A_KE_K(\xi,\eta)-[B_K(\xi,\eta)]^2)}}\right)\right]\right\}\\[0.5ex]
\end{multlined}
\label{eqn:55}
\end{equation}
Case 2:
In this case $(\xi,\eta)$ lies on the boundary element but doesn't lie on the boundary on which we are integrating. Hence,
\begin{equation}
\begin{multlined}
F_1^K(\xi,\eta)\\
=\frac{l^K}{2\pi}\left\{ln(l^K)+\left(1+\frac{B_K(\xi,\eta)}{2A_K}\right)ln\left|\frac{B_K(\xi,\eta)}{2A_K}\right|\right.\\
\left.-\frac{B_K(\xi,\eta)}{2A_K}ln\left|\frac{B_K(\xi,\eta)}{2A_K}\right|-1\right\}
\end{multlined} 
\label{eqn:56}
\end{equation}
$F_2^K(\xi,\eta)$ - case 1:\\
In this case $(\xi,\eta)$ doesn't lie on the boundary element on which we are performing integration. So integrating analytically, it can be stated that 
\begin{equation}
\begin{multlined}
F_2^K(\xi,\eta)=\frac{l^k[n_x^K(x^K-\xi)+n_y^K(y-\eta)]}{\pi\sqrt{(4A_KE_K(\xi,\eta)-[B_K(\xi,\eta)]^2)}}\\[1ex]\times\left[tan^{-1}\left(\frac{2A_K+B_K(\xi,\eta)}{\sqrt{(4A_KE_K(\xi,\eta)-[B_K(\xi,\eta)]^2)}}\right)\right.\\[1ex]
\left.-tan^{-1}\left(\frac{B_K(\xi,\eta)}{\sqrt{(4A_KE_K(\xi,\eta)-[B_K(\xi,\eta)]^2)}}\right)\right]
\end{multlined}
\label{eqn:57}
\end{equation}
Case 2:\\
In this case $(\xi,\eta)$ lies on the boundary element but doesn't lie on the boundary on which we are integrating. Hence
\begin{equation}
F_2^K(\xi,\eta)=0
\label{eqn:57a}
\end{equation}
Based on multiple parameters of the system, which has been detailed in section II, a one-time FEM on the entire domain via "COMSOL MULTIPHYSICS" software has been performed. Using the FEM, magnetic vector potential values on both circles, i.e. at $(\bar{x}^K,\bar{y}^K)$, k=1, 2, ......, 80 have been obtained. The output can be represented as 

\begin{align*}
\begin{pmatrix}
\bar{A}_1 \\\vdots\\\bar{A}_{80}
\end{pmatrix}_{FEM}
\end{align*}
The values obtained through COMSOL is exported to the MATLAB environment. Using the FEM-BEM coupling condition, the output of FEM has been used as the input of BEM.
\begin{align*}
\begin{pmatrix}
\bar{A}_1 \\\vdots\\\bar{A}_{80}
\end{pmatrix}_{BEM}
=
\begin{pmatrix}
\bar{A}_1 \\\vdots\\\bar{A}_{80}
\end{pmatrix}_{FEM}
\end{align*}

Then, substituting (\ref{eqn:49}) in (\ref{eqn:48}) in $\bar{P}_K$ for k=1, 2, 3, ........., 80 have been obtained. After getting this complementary solution set, $A(\xi,\eta)$ at any point in the interior part of the annulus can be calculated by (\ref{eqn:59a})

\begin{equation}
\begin{multlined}
A(\xi,\eta)=
\begin{pmatrix}
\bar{A}_1 &\cdots&\bar{A}_{80}
\end{pmatrix}
\begin{pmatrix}
F_2^1(\xi,\eta) \\\vdots\\F_2^{80}(\xi,\eta)
\end{pmatrix}\\[1ex]
-
\begin{pmatrix}
\bar{P}_1 &\cdots&\bar{P}_{80}
\end{pmatrix}
\begin{pmatrix}
F_1^1(\xi,\eta) \\\vdots\\F_1^{80}(\xi,\eta)
\end{pmatrix}
\end{multlined}
\label{eqn:59}
\end{equation}
In a more abridged form, (\ref{eqn:59}) can be written as 
\begin{equation}
A(\xi,\eta)=\sum_{K=1}^{80}\left\{\bar{A}_KF_2^{K}(\xi,\eta)-\bar{P}_KF_1^{K}(\xi,\eta)\right\}
\label{eqn:59a}
\end{equation}
Where $(\xi,\eta)$ $\in$ R. Hence, in further calculating the magnetic vector potential at any interior point of the domain, the computation process will solely depend upon BEM and will be free from the adaptive domain discretization process. In the modified approach, subsequent coding has been developed in a MATLAB environment to calculate magnetic vector potential at different points of the aluminum plate.

\section{ LEVITATING SYSTEM DESCRIPTION}
\begin{figure}[h]
\includegraphics[width=4cm]{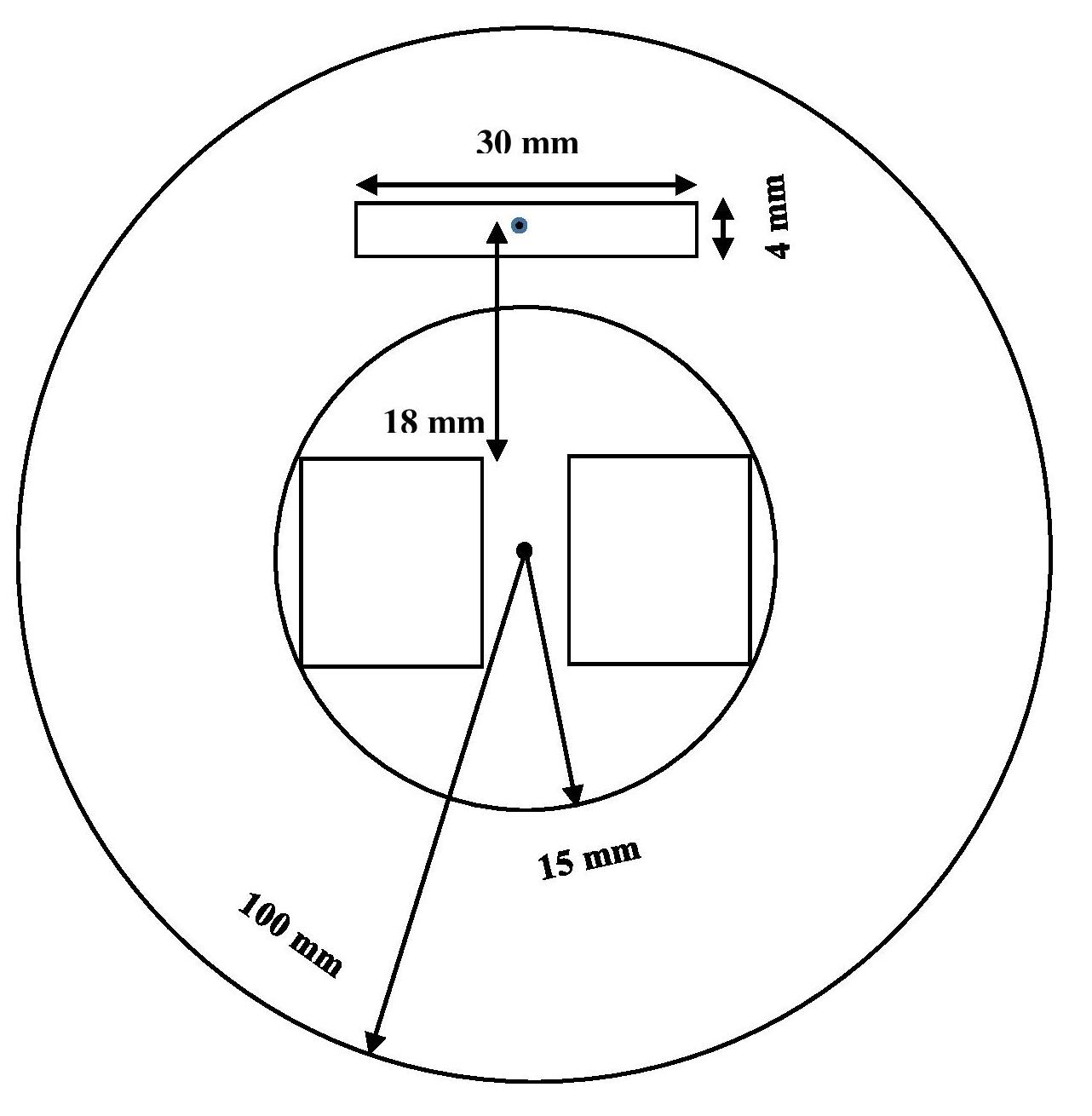}
\centering
\caption{ Magnetic levitation set up at t = 0\cite{jena}}
\label{fig:set1}
\end{figure}

\begin{figure}[h]
\includegraphics[width=4cm]{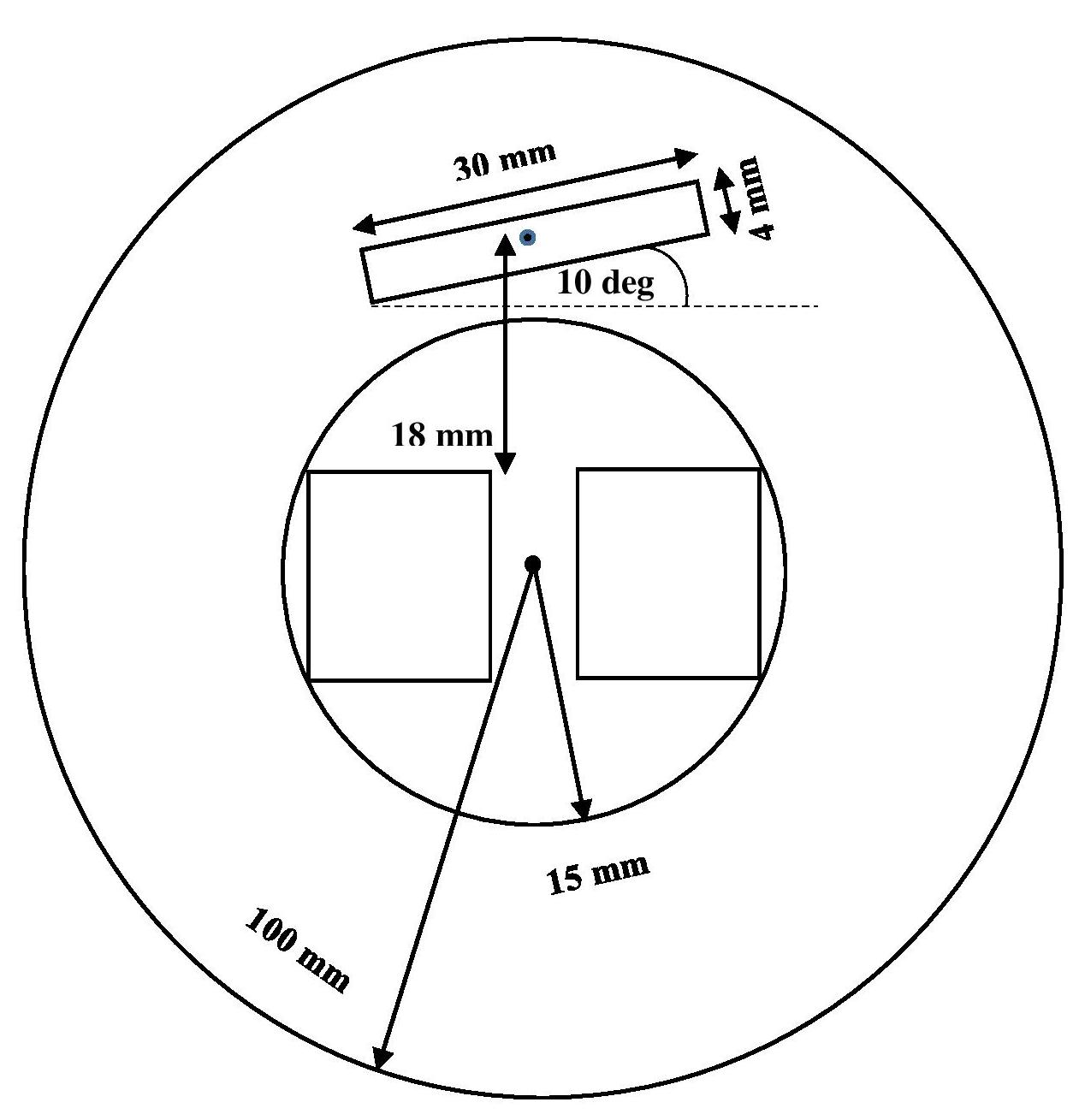}
\centering
\caption{Disturbed levitation set up after some instant\cite{jena}}
\label{fig:set2}
\end{figure}

The simulation setup of the problem follows Fig. \ref{fig:set1}. A cuboidal plate of aluminum of conductivity ($\sigma$=3.7$\times$$10^7$$\mho$m), having a depth of 1 mm, is placed above one cylindrical coil (incoming and outgoing) having 960 turns. Coils are shown by domains that carry a homogeneous azimuthal current density for numerical purposes. The coils are made of copper wire with 1.2 mm diameter and contain insulating layers. The levitation height 'Y' denotes to the distance between the midpoint of the lower edge of the plate and the upper edge of the current-carrying area (Y = 0) for t $\leq$ 0 the plate lies above the coils at a distance of 18 mm due to the thickness of the winding form. Centre of gravity of the system at (x = 0 mm, y = 18 mm). For calculation purposes, two fictitious circles are defined; the inner one has a radius of 15 mm, and the outer ring has a radius of 100 mm. After some instant due to oscillations for the applied force, the plate shifts 2 mm upward and 2 mm left-wise and rotated at an angle of 10 (refer to final position) as shown in Fig. \ref{fig:set2} (not to scale).For t $\geq$ 0, sinusoidal current\emph{ i(t)} flows in the coil segment in opposite directions.

\begin{equation}
i(t) = \hat{i}\,sin(2\pi ft), \hat{i}=200A, f_0=50Hz
\label{eqn:it}
\end{equation}

The interaction of the induced eddy current on current transferring conductor is omitted for easy calculation purpose. For this reason, the current in (\ref{eqn:it}) can be regarded as impressed.

\section{RESULTS AND DISCUSSIONS}
\begin{figure}[H]
\includegraphics[width=8cm]{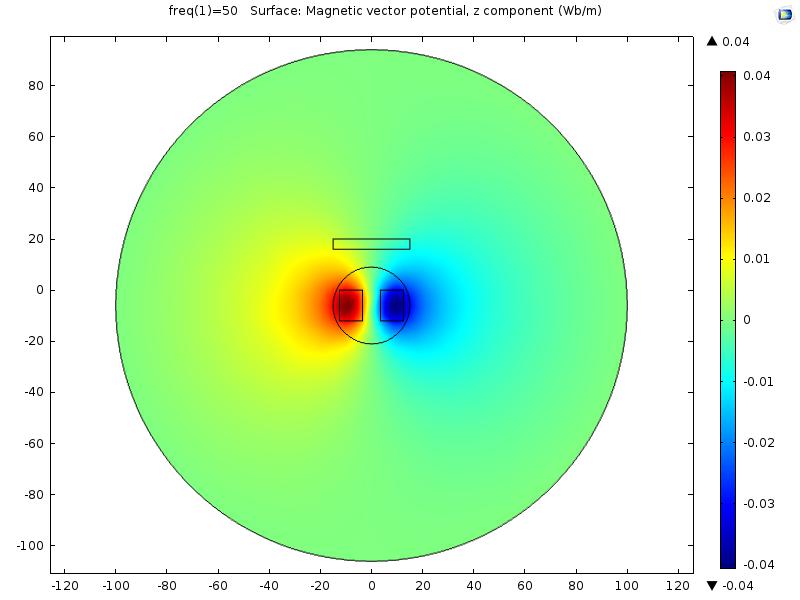}
\centering
\caption{Magnetic vector potential distribution at initial condition \cite{jena}}
\label{fig:sd}
\end{figure}

\begin{figure}[H]
\includegraphics[width=8cm]{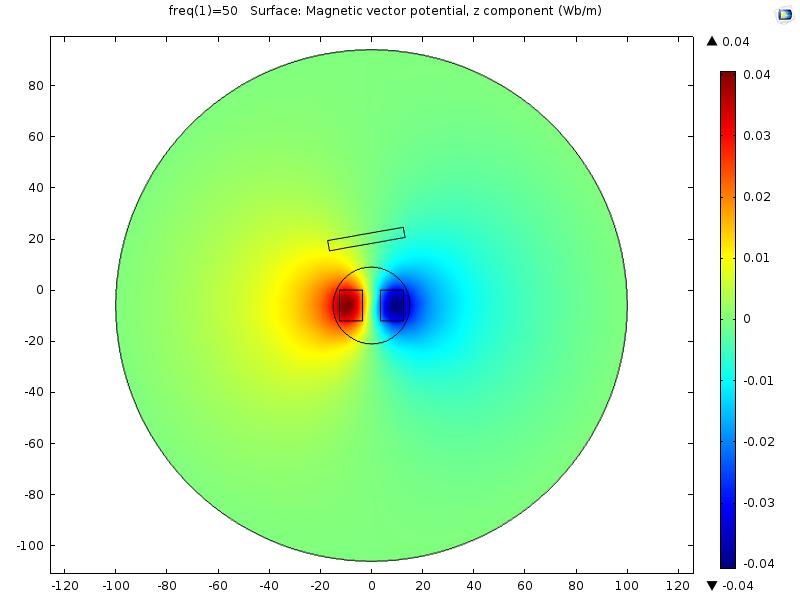}
\centering
\caption{Magnetic vector potential distribution at some instant after applying disturbance \cite{jena}}
\label{fig:dd}
\end{figure}

\begin{table}[H]
\caption{INITIAL POSITION (40 NODES)}
\centering
\begin{tabular}{c rrr}
\hline\hline
Measured  &Calculated  &Error&Average \\
value (T m)&value (T m) &&error\\[0.5ex]
\hline\\[0.2ex]
-0.0065 & -0.0066 & 1.53\%&7.47\%\\[0.5ex]
-0.0025&-0.0028&12\%&\\[0.5ex]
0.0024&0.0025&4.16\%&\\[0.5ex]
0.0063&0.0062&1.58\%&\\[0.5ex]
-0.0090&-0.0085&5.55\%&\\[0.5ex]
-0.0046&-0.0040&13.04\%&\\[0.5ex]
0.0042&0.0037&11.90\%&\\[0.5ex]
0.0090&0.0081&10\%&\\[0.5ex]
\hline 
\end{tabular}
\label{tab:table1}
\end{table}
 
\begin{table}[H]
\caption{DISTURBED POSITION (40 NODES)}
\centering
\begin{tabular}{c rrr}
\hline\hline
Measured  &Calculated  &Error&Average \\
value (T m) &value (T m) &&error\\[0.5ex]
\hline\\[0.2ex]
-0.0045 & -0.0045 & 0\%&5.87\%\\[0.5ex]
-0.0013&-0.0013&0\%&\\[0.5ex]
0.0026&0.0030&15.38\%&\\[0.5ex]
0.0067&0.0069&2.98\%&\\[0.5ex]
-0.0067&-0.0063&5.97\%&\\[0.5ex]
-0.0031&-0.0025&19.35\%&\\[0.5ex]
0.0044&0.0044&0\%&\\[0.5ex]
0.0092&0.0089&3.26\%&\\[1ex]
\hline 
\end{tabular}
\label{tab:table2}
\end{table}

\begin{figure}[H]
\includegraphics[width=6cm]{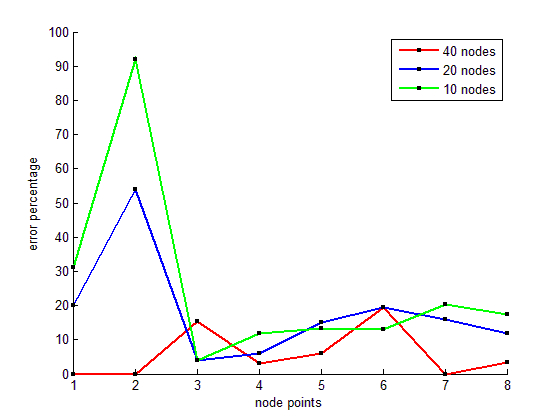}
\centering
\caption{Local error analysis at initial time\cite{jena}}
\label{fig:locs}
\end{figure}

\begin{figure}[H]
\includegraphics[width=6cm]{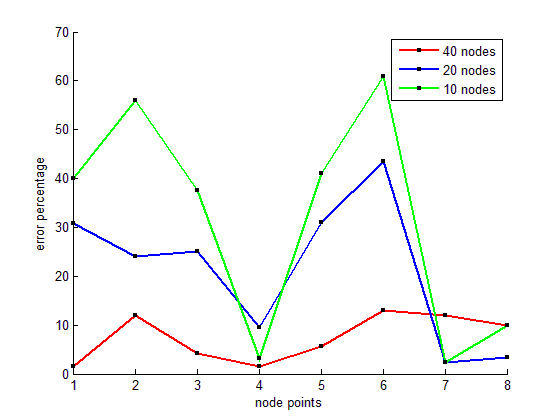}
\centering
\caption{Local error analysis after displacement \cite{jena}}
\label{fig:locd}
\end{figure}

\begin{figure}[H]
\includegraphics[width=6cm]{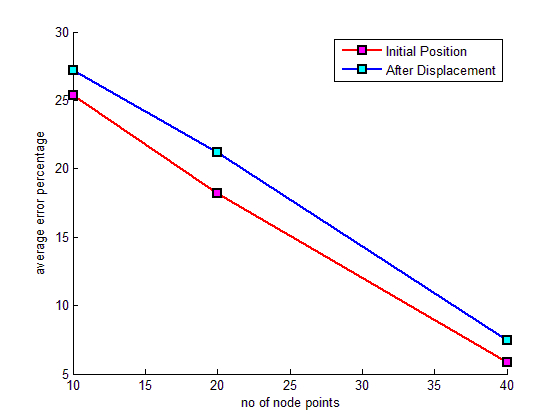}
\centering
\caption{Relative analysis of average error with number of nodes \cite{jena}}
\label{fig:gross}
\end{figure}

A comparative analysis of the magnetic vector potential obtained from Multi-physics software COMSOL using re-meshing strategies as per Fig.\ref{fig:sd}-\ref{fig:dd} and computed results from MATLAB are shown in Table \ref{tab:table1}-\ref{tab:table2} obtained with the help of developed BEM-FEM code. Relative error analysis graphs as in Fig. \ref{fig:locs}-\ref{fig:locd} are also attached to know the behavior of error with a variation of nodal points. With an increasing number of node points on the boundary of the impressed current domain error, it can be seen that error is decreased at a sharp rate. This work is evaluated in Dell Inspiron 15 series computer operating with 2.1 GHz processor with 4 GB RAM. The interface between the aluminum plate and air coincides with coupling boundary $\Gamma_{BEM}$=$\Gamma_{BEM}$ and it has 8 nodes. The evaluated magnetic vector potential on the boundary of the aluminum body serves as input for tangential and normal force density and subsequently force and torque calculation for Maxwell stress tensor-based method during comparison of hardware and software-based result. The drawback of the work is that we assumed a 'constant approximation' in gradient calculation. If we assume 'linear' or 'quadratic approximation', the error between the measured values and calculated values will be very less. A fictitious circle is drawn around a metal plate and current-carrying conductor domain. However, it may not be possible to draw a circle in between them due to constrain in intermediate distance. 

\section{CONCLUSION}
Magnetic vector potential on the boundary of the conducting surface has been calculated using the modified FEM-BEM coupling approach. A  good agreement between the calculated results from two different software has been observed. The efficiency of the method has been illustrated by comparing the results obtained from nodal performance analysis. Previously, coding of BEM-FEM coupling approach has been done in FORTRAN environment, MATLAB based coding is newly introduced in this paper. Significant savings in computation time and storage requirements are achieved. Magnetic forces, normal force density, and tangential force density can be calculated from magnetic vector potential as per the standard formula from the value obtained on the boundary. This concept is extended to compute net force and torque by simply integrating force and torque over the boundary, assuming that the boundary is a coupled one. Then only the analysis of the dynamic behavior of the levitating body is possible.

\ifCLASSOPTIONcaptionsoff
  \newpage
\fi



\begin{thebibliography}{1}

\bibitem{karl}
H.~Karl, J.~Fetzer, S.~Kurz, G.~Lehner and W.~M. Rucker, ``Description 
      of TEAM Workshop Problem 28: An Electrodynamic Levitation device''.

\bibitem{kurz1}
S.~Kurz, J.~Fetzer and G.~Lehner, ``Three-dimensional 
      transient BEM-FEM coupled analysis of electrodynamic levitation 
      problems,'' \textit{IEEE Trans. Magn.}, vol. 32, no. 3, pp. 1062-1065, 1996.

\bibitem{fetter}
J.~Fetter \emph{et al}, ``Transient BEM-FEM coupled analysis of 3-D electromechanical systems: a watch stepping motor driven by a thin wire coil," \textit{IEEE Trans. Magn.}, vol. 34, no. 5, pp. 3154-3157, 1998.

\bibitem{suuriniemi}
S. Suuriniemi, K. Forsman, L. Kettunen and J. Makinen, "Computation of eddy currents coupled with motion," \textit{IEEE Trans. Magn.}, vol. 36, no. 4, pp. 1341-1345, July 2000.


\bibitem{jena}
A.~Jena, S.~Sarkar and Mashuq-un-Nabi, "A modified BEM-FEM coupling approach for 3D electromagnetic levitation problem," 2015 Annual IEEE India Conference (INDICON), New Delhi, 2015, pp. 1-6.

\bibitem{nicolet} 
A.~Nicolet, F.~Delinc\'e, A.~Genon and W.~Legros, ``Finite elements-
      boundary elements coupling for the movement modelling in two-    
      dimensional structures," \textit{J. Phys. III}, vol. 2, pp. 2035-2044, 1992.

\bibitem{kurz2}
S.~Kurz, J.~Fetzer, G.~Lehenr, W.~M. Rucker, ``A novel formulation for 3D eddy current problems 
      with moving bodies using a Lagrangian description and BEM-FEM coupling," \textit{IEEE Trans. Magn.}, vol. 34, no. 5,  pp. 3068-3073, 1998.

\bibitem{fetzer}
J.~Fetzer, S.~Kurz, G.~Lehner and W.~M. Rucker, ``Analysis of an actuator with eddy currents and iron saturation: comparison between a FEM and a BEM-FEM coupling approach," \textit{IEEE Trans. Magn.}, vol. 35, no. 3, pp. 1793-1796, 1999.

\bibitem{kurz3}
S.~Kurz, J.~Fetzer, G.~Lehner and W.~M. Rucker, ``The application of the BEM-FEM coupling method for the treatment of three-dimensional nonlinear shielding problems of low-frequency fields using the example of the TEAM problem 21," \textit{Archive for Electrical Engineering}, vol. 80, no. 2, pp. 91-104, 1997.

\bibitem{ang}
W.~T. Ang, A Beginner’s Course in Boundary Element Method, Universal Publisher, 2007. 





\end{thebibliography}
%




\end{document}